\documentclass[a4paper,11.5pt,reqno]{amsart}

\usepackage{graphicx}
\usepackage[utf8]{inputenc}
\usepackage{amsmath,amssymb,amsthm}
\usepackage{enumitem}
\usepackage{xcolor}
\usepackage{hyperref}

\usepackage{lineno}
\modulolinenumbers[5]
\numberwithin{equation}{section}

\hypersetup{colorlinks = true, linkcolor = blue, urlcolor = blue, citecolor = blue}

\theoremstyle{plain}

\theoremstyle{definition}

\theoremstyle{remark}

\newcommand{\oeis}[1] { \href {https://oeis.org/#1} {\path{#1}} }

\newcommand{\Rx} { \mathbb{R}^\times }

\begin{document}

\title[Somos-4 and a quartic Surface in $\mathbb{RP}^{3}$]
      {Somos-4 and a quartic Surface in $\mathbb{RP}^{3}$}

\author[H. Ruhland]{Helmut Ruhland}
\address{Santa F\'{e}, La Habana, Cuba}
\email{helmut.ruhland50@web.de}



\begin{abstract}
The Somos-4 equation defines the sequences with this name. Looking at these sequences with an
additional property we get a quartic polynomial in 4 variables. This polynomial defines a rational, projective surface in $\mathbb{RP}^{3}$. Here some generators of the subgroup of $Cr_3 (\mathbb{R})$ are determined, whose birational maps are automorphisms of
the quartic  surface.
\end{abstract}

\date{\today}

\maketitle

\section{Introduction}

In this article I define a rational, quartic surface in $\mathbb{RP}^{3}$ related to a special Somos-4 sequence.
The special property is: the 2 subsequences with even and odd indices are also Somos-4 sequences.
In the following to be short I write \emph{"even and odd subsequences"} instead of the ... above.
The birational maps in $Cr_3 (\mathbb{R})$ generating the automorphism group
of this variety are given. \\

Here an outline of this article: \\
In section \ref{SomosSequences} the Somos-4 sequences are defined and the here relevant properties are given. The transformation $T$ transforming Somos sequences into Somos sequences is defined.
In the next section \ref{SomosReps} a representation of the transformation group $T$ in $Cr_{n-1} (\mathbb{R})$ for general Somos-$n$ sequences is given. \\
In section \ref{SomosSequencesEvenOdd} the condition for the even and odd subsequences to be Somos-4 sequences is derived and it defines a rational, quartic surface. In the given example with only integer entries the odd subsequence is the classical Somos-4 sequence with initial values $1, 1, 1, 1$.
From this example follows a divisibility relation between this classical Somos-4 sequence and \emph{another} Somos-4 sequence.
In section \ref{Sym} birational automorphisms of this surface are determined. \\

As appendices the representations of the group of transformations of Somos sequences in the Cremona
groups are given for Somos-2, Somos-3 and Somos-5. In the case of Somos-3 and the representation in the Cremona group of rank 2 the invariant curves and the birational automorphisms are determined in appendix \ref{RepSomos3}.

\section{About Somos-4 sequences and the group $T$ of transformations \label{SomosSequences}}

A Somos-4 sequence $\dots, a_{-1}, a_0, a_1, a_2, a_3, a_4, a_5, \dots$ with the indices $n \in \mathbb{Z}$ fulfills the Somos-4 equation:
\begin{equation}
	a_{n} a_{n+4} = a_{n+1} a_{n+3} + a_{n+2}^2  \label{Somos4Equ}
\end{equation}

There are 3 transformations, which transform general Somos sequences (not only Somos-4 sequences) into Somos sequences, due to a common property of all Somos equations. Let $Seqs$ be the set of all Somos sequences. The transformations are $M(b,c) : Seqs \rightarrow Seqs, a_{n} \mapsto a_{n} \, b \, c^n$, the reflection transform $R : Seqs \rightarrow Seqs, a_{n} \mapsto a_{-n}$ and the shift transform $F : Seqs \rightarrow Seqs, a_{n} \mapsto a_{n + 1}$.
$R$ and $F$ generate an infinite dihedral group $D_\infty$. \\

Let "$\times$" denote the direct product of groups, "$\rtimes$" denote the semidirect product of groups
i.e. split group extensions. "$\rtimes$" has higher priority than "$\times$", i.e $\mathbb{Z} \times \mathbb{Z} \rtimes \mathbb{Z}_2$ without parentheses means $\mathbb{Z} \times (\mathbb{Z} \rtimes \mathbb{Z}_2)$. $\Rx$ is the multiplicative group in $\mathbb{R} \setminus \{ 0 \}$. \\

The group of transforms of Somos sequences is then; 
\begin{equation}
	T = \langle \, M(b,c), F, R \, \rangle \simeq ( \, (\Rx \times \Rx) \rtimes \mathbb{Z} \, ) \rtimes \mathbb{Z}_2  \label{GroupSomosTrans}
\end{equation}

Because $\langle \, M(b,c) \, \rangle \simeq \Rx \times \Rx$ and $\langle  \, F \, \rangle \simeq \mathbb{Z}$ do not commute, the commutator $[M(b,c), F] \in \langle \, M(b,1) \, \rangle$, the group $\langle \, M(b,c), F \, \rangle$ above is the semidirect product $(\Rx \times \Rx) \rtimes \mathbb{Z}$. \\  
$T$ has normal subgroups: the center $\langle \, M(b,1) \, \rangle$, $\langle \, M(b,1), F^n \, \rangle$
\dots. These normal subgroups define the quotient groups:
\begin{align}
	T^* & = T \setminus \langle \, M(b,1) \, \rangle 
           \simeq (\Rx \times \mathbb{Z}) \rtimes \mathbb{Z}_2    \label{GroupSomosTransStar} \\
	T^{*}_n & = T \setminus \langle \, M(b,1), F^n \, \rangle 
           \simeq (\Rx \times \mathbb{Z}_n) \rtimes \mathbb{Z}_2  \label{GroupSomosTransStar_n}
\end{align}
\begin{align}
	T^{*}_1 \simeq \Rx \rtimes \mathbb{Z}_2 \qquad  T^{*}_2 \simeq \mathbb{Z}_2 \times \Rx \rtimes \mathbb{Z}_2  \label{GroupSomosTransStar_12}
\end{align}
The later three $^{*}$-groups have the following infinite subgroups: \\
The group $T^{*}$ has the \emph{discrete} infinte dihedral $\mathbb{Z} \rtimes \mathbb{Z}_2$ and another \emph{continuous} infinte dihedral $\Rx \rtimes \mathbb{Z}_2$ as subgroups. $T^{*}_2$ and $T^{*}_1$ have
only the continuous infinte dihedral subgroup as infinte subgroup. \\

Because the Somos-4 equation \ref{Somos4Equ} is linear in the terms with the highest and lowest indices $a_{}, a_{n+4}$, we get the following 2 rational recurrences:
\begin{equation}
	a_{n} = (a_{n-1} a_{n-3} + a_{n-2}^2) / a_{n-4} \qquad 
   a_{n} = (a_{n+1} a_{n+3} + a_{n+2}^2) / a_{n+4} \label{SomosRec}
\end{equation}
Given 4 subsequent terms of a sequence, these 2 recurrences allow us to calculate all other terms with higher and lower indices. Let our 4 subsequent initial terms be $a_0, a_1, a_2, a_3$.
With the first recurrence $a_4 = (a_1 a_3 + a_{2}^2) / a_0, a_5 = (a_0 a_3^2 + a_1 a_2 a_3 + a_2^3) / (a_0 a_1), \dots$ can be determined. With the second recurrence $a_{-1} = (a_0 a_2 + a_{1}^2) / a_3, a_{-2} = (a_0^2 a_3 + a_0 a_1 a_2 + a_1^3) / (a_2 a_3), \dots$ can be determined. For Somos-3, ..., Somos-7 sequences the denominators of all terms $a_n$ are monomials in the initial values, the so called Laurent property, see \cite{FoZe}. \\

For a history of Somos sequences, see e.g. \cite{Gale}. \\

\section{Somos-$n$: A representation of $T$ in the Cremona group $Cr_{n-1} (\mathbb{R})$ \label{SomosReps}}

In the previous section it was shown for Somos-4 sequences, that all sequence terms can be expressed
by the 4 initial values  $a_0, a_1, a_2, a_3$ as a homogenous, rational expression with the Laurent property. In general Somos-$n$ sequences all terms $a_i$ can be expressed by $n$ initial values.
The Laurent property is fulfilled only for $2 \le n \le 7$. \\

Let $t \in T$ be a transformation of Somos-$n$ sequences. With these rational expressions for the sequence terms we can define a representation of $T$ in the Cremona group $Cr_{n-1} (\mathbb{R})$, $Rep : T \rightarrow Cr_{n-1} (\mathbb{R})$. The image of $t$ is:
\begin{align} \begin{split}
Rep (t) : \enspace  \mathbb{RP}^{n-1} \dashrightarrow \enspace \mathbb{RP}^{n-1}, 
 \; & ( a_0 \, : a_1 \, : \dots \, : a_{n-2} \, : a_{n-1} ) \\
 \mapsto & \enspace  ( t (a_0) \, : t (a_1) \, : \dots \, : t (a_{n-2}) \, : t (a_{n-1}) ) \label{tImageRep}
\end{split} \end{align} \\

The representation of $T$ is not faithful, because $T$ has a non-trivial center and the image of $T$ acts
on a projective space. But for $n \ge 4$ the representation of $T^{*} = T \; mod$ its center is faithful, see appendix \ref{RepSomos5}. For $n = 3$ the representation of $T^{*}_2$ is faithful,  see appendix \ref{RepSomos3}. For $n = 2$ the representation of $T^{*}_1$ is faithful, see appendix \ref{RepSomos2}.

\section{Even and odd Somos-4 subsequences \label{SomosSequencesEvenOdd}}

Looking for sequences, whose 2 even and odd subsequences are also Somos-4, we get this quartic polynomial
in 4 subsequent entries as condition:
\begin{equation}
	S_{n} = a_{n+0}^2 a_{n+3}^2 + a_{n+1}^2 a_{n+2}^2   + a_{n+0} a_{n+2}^3 + a_{n+3} a_{n+1}^3 
          + 2 \, a_{n+0} a_{n+1} a_{n+2} a_{n+3} \label{Cond}
\end{equation}
This polynomial occurs as a factor in $a_{n} a_{n+8} - a_{n+2} a_{n+6} - a_{n+4}^2$. Before factoring
the $a_{n+4}, a_{n+6}, a_{n+8}$ are expressed by the four subsequent $a_{n}, a_{n+1}, a_{n+2}, a_{n+3}$.
Factoring $S_{n+1}$ (expressing $a_{n+4}$ by the four subsequent $a_{n}, a_{n+1}, a_{n+2}, a_{n+3}$)
we get this equation:
\begin{equation}
	S_{n+1} = S_{n} (a_{n+1} a_{n+3} + a_{n+2}^2) / a_{n+0}^2 \label{zzz}
\end{equation}
If $S_0 = 0$ all following $S_1, S_2, \dots$ are $0$ and the even and odd subsequences are Somos-4 sequences. \\

An example for a Somos-4 sequence with this extra property with \emph{only integers} as entries is
\oeis {A006769} in N.J.A. Sloane's On-Line Encyclopedia of Integer Sequences, OEIS\textregistered. \\

All four subsequent entries of this sequence A006769 define the integer coordinates of points
on the quartic $S_0$, so we have an infinite number of points with integer coordinates on the surface:
i.e. $( 0 \, : 1 \, : 1 \, : -1 ), ( 1 \, : 1 \, : -1 \, : 1 ), \dots, $
$( 7 \, : -4 \, : -23 \, : 29 ), \dots$. \\

The even subsequence is \oeis {A051138}, the odd subsequence is the classical (initial values $1, 1, 1, 1$) Somos-4 \oeis {A006720}, starting with index $n = 2$ and every entry multiplied by $(-1)^{n}$. \\

Divisibility in sequences: \\
Because \oeis {A006769} is a strong (elliptic) divisibility sequence i.e.
\oeis {A006769}$\, (n) \; \vert$ \\ \oeis {A006769}$\, (n k)$, this divisibilty relation
induces divisibilty relations between subsequences with indices in an arithmetic progression and so divisibility in the 2 even and odd subsequences:
\begin{itemize}
   \item \oeis {A006720}$\, (n) \; \vert$ \oeis {A006720}$\, (n + (2n - 3) k)$, given already in
         a comment in A006720 by Peter H. van der Kamp, 2015.
   \item \oeis {A006720}$\, (n) \; \vert$ \oeis {A051138}$\, ((2n - 3) k)$, now between
         2 \emph{different} Somos-4 sequences.
\end{itemize}

Question: Do there exist other (than in the example above given) Somos-4 sequences with integer entries
          and even and odd Somos-4 subsequences?

\section{The quartic surface and its birational symmetry group \label{Sym}}

Now we take $S_0$ as polynomial $S$ defining the quartic surface: 
\begin{equation}
	S = S_{0} = a_{0}^2 a_{3}^2 + a_{1}^2 a_{2}^2   + a_{0} a_{2}^3 + a_{3} a_{1}^3
          + 2 \, a_{0} a_{1} a_{2} a_{3} \label{QuartSurface}
\end{equation}

The following hint is from Igor Dolgachev: \\
The quartic surface $S$ is rational, so its group of birational automorphisms coincides with the whole Cremona group of rank 2. To see that it is rational, project the surface from its singular point
$( 1 \, : 0 \, : 0 \, : 0 )$. The surface becomes the double cover of the plane branched along a curve of degree 6 with an ordinary singular point of multiplicity $5$. The surface is not rational if and only if the singular points of the branch sextic are ADE rational double points. \\

The 3 transformations in section \ref{SomosSequences} now appear as birational symmetry maps
in $Cr_3 (\mathbb{R})$ leaving $S$ invariant. \\

The transform $a_{n} \mapsto a_{n} \, b c^n$ results in a 1-paramter map depending on $c \ne 0$:
\begin{equation}
M(c) : \enspace  \mathbb{RP}^3 \dashrightarrow \enspace \mathbb{RP}^3, 
 ( a_0 \, : a_1 \, : a_2 \, : a_3 ) \mapsto \enspace ( a_0 \, : a_1 c \, : a_2 c^2  \, : a_3 c^3 ) \label{MapMonomial}
\end{equation}

The reflection transform $a_{n} \mapsto a_{3-n}$ results in the reflection map $R$:
\begin{equation}
R : \enspace  \mathbb{RP}^3 \dashrightarrow \enspace \mathbb{RP}^3, 
 ( a_0 \, : a_1 \, : a_2 \, : a_3 ) \mapsto \enspace ( a_3 \, : a_2 \, : a_1 \, : a_0 ) \label{MapRefl}
\end{equation}

The shift transform $a_{n} \mapsto a_{n+1}$ results in the shift map $F$:
\begin{equation}
F : \enspace  \mathbb{RP}^3 \dashrightarrow \enspace \mathbb{RP}^3, 
 ( a_0 \, : a_1 \, : a_2 \, : a_3 ) \mapsto \enspace ( a_1 \, : a_2 \, : a_3  \, : (a_1 a_3 + a_{2}^2) / a_0 ) \label{MapShift}
\end{equation}

The group $\langle \, R, F \, \rangle$ is isomorph to the infinite dihedral group $D_\infty$.
This group can also be generated
by $R$ and an additional involution $G = R F$, so $\langle \, R, F \, \rangle = \langle \, R, G \, \rangle$: 
\begin{equation}
G = R F : \enspace  \mathbb{RP}^3 \dashrightarrow \enspace \mathbb{RP}^3, 
 ( a_0 \, : a_1 \, : a_2 \, : a_3 ) \mapsto \enspace ( (a_1 a_3 + a_{2}^2) / a_0 \, : a_3 \, : a_2  \, : a_1 ) \label{Invo}
\end{equation}
This map stems from the reflection transformation $a_{n} \mapsto a_{4-n}$. As polynomial map it is of degree $2$. The determinant of the Jacobian $\det\mathrm{jac}$ is the union of the plane $a_0 = 0$ and the quadratic surface $a_1 a_3 + a_{2}^2 = 0$. \\

The 1-parameter map $M(c)$ does not commute with the 2 involutions $R, G$. We have $R M(c) R^{-1} = G M(c) G^{-1} = M(c)^{-1} = M(c^{-1})$. So $\langle \, M(c), R \, \rangle$ and $\langle \, M(c), G \, \rangle$ are infinite, \emph{now continuous} dihedral groups $D_\infty (c)$. \\

$\hat{T} = \langle \, M(c), \; R, G \, \rangle$ defines a representation of the group of Somos sequence transformations $T$, see \ref{GroupSomosTrans}, in $Cr_3 (\mathbb{R})$. This
representation is not faithful and $\hat{T}$ is isomorph to quotient of $T$, $\hat{T} \simeq T^* \simeq (\mathbb{R} \times \mathbb{Z}) \rtimes \mathbb{Z}_2$, see \ref{GroupSomosTransStar}.
The degrees of the maps $F^n, M(c), R F^n$ and all group elements are $1, 2, 3, 5, 8, 10, 14, 18, \dots$
for $n \ge 0$.  \\

Question: Are there besides the quartic surface $S$ other surfaces invariant under $\hat{T}$ or a nontrivial subgroup like $\mathbb{Z} \rtimes \mathbb{Z}_2$, the infinte dihedral group? \\

Because $S$ is quadratic in $a_0$ and $a_3$ we can construct two further birational involutions.
This is done in a similar manner as its done e.g. for the groups acting on the Markov triples and
on the curvatures in an Apollonian circle packing.
\begin{equation}
H : \enspace  \mathbb{RP}^3 \dashrightarrow \enspace \mathbb{RP}^3, 
 ( a_0 \, : a_1 \, : a_2 \, : a_3 ) \mapsto \enspace ( (a_2^3 - 2 a_1 a_2 a_3) / a_3^2 - a_0 \, : a_1 \, : a_2  \, : a_3 ) \label{MapExtra}
\end{equation}
This map as polynomial map is of degree $3$ and commutes with the 1-parameter map $M(c)$.
$\det\mathrm{jac}$ is the plane $a_3 = 0$ with multiplicity $8$. The other involution belonging to $a_3$ is just this $H$ conjugated by $R$. \\

Now the symmetry group of $S$ is $\langle \, M(c), \; R, G, H \, \rangle = \langle \, \hat{T}, H \, \rangle \subset Cr_3 (\mathbb{R})$. \\

\vspace{1.0cm}
\newpage

\noindent \textbf{\large Appendices}

\appendix

\section{Somos-2: The representation of $T$ in $Cr_1 (\mathbb{R})$ \label{RepSomos2}}

The terms of the Somos-2 sequence can be expressed in a simple way by the 2 initial values $a_0, a_1$ as  $a_n = a_1^{n} / a_0 ^{n-1}$ for all $n$. \\  

The transform $a_{n} \mapsto a_{n} \, b c^n$ results in a 1-paramter map depending on $c \ne 0$:
\begin{equation}
M(c) : \enspace  \mathbb{RP}^1 \dashrightarrow \enspace \mathbb{RP}^1, 
 ( a_0 \, : a_1 ) \mapsto \enspace ( a_0 \, : a_1 c ) \label{MapMonomialP1}
\end{equation}

The reflection transform $a_{n} \mapsto a_{1-n}$ results in the reflection map $R$:
\begin{equation}
R : \enspace  \mathbb{RP}^1 \dashrightarrow \enspace \mathbb{RP}^1, 
 ( a_0 \, : a_1 ) \mapsto \enspace ( a_1 \, : a_0 ) \label{MapReflP1}
\end{equation}

The shift transform $a_{n} \mapsto a_{n+1}$ results in the shift map $F$:
\begin{equation}
F : \enspace  \mathbb{RP}^1 \dashrightarrow \enspace \mathbb{RP}^1, 
 ( a_0 \, : a_1 ) \mapsto \enspace ( a_1 \, : a_1^2 / a_0 ) = ( a_0 \, : a_1 ) a_1 / a_0 \label{MapShiftP1}
\end{equation}
Because we work in a projective space $\mathbb{RP}$, a common factor on a vector is an equivalenve. In this case F is \emph{trivial and the unit element} of order $1$. \\

The group $\langle \, R, F \, \rangle$ is therefore isomorph to the finite dihedral group $D_2$. \\

The 1-parameter map $M(c)$ does not commute with the involution $R$. We have $R M(c) R^{-1} = M(c)^{-1} = M(c^{-1})$. \\

$\hat{T} = \langle \, M(c), \; R \, \rangle \simeq T^{*}_1 \simeq \Rx \rtimes \mathbb{Z}_2$, see \ref{GroupSomosTransStar_12}, defines a representation of the group of Somos sequence transformations $T$ in $Cr_1 (\mathbb{R})$. \\

The generators $M(c), R$ as linear fractional maps in $\mathbb{R}$ with $x = a_0 / a_1$ are:
\begin{align}
M(c) & : \enspace  \mathbb{R} \rightarrow \enspace \mathbb{R}, 
 x \mapsto \enspace x \, c \label{MapMonomialP1LF} \\
R & : \enspace  \mathbb{R} \rightarrow \enspace \mathbb{R}, 
 x \mapsto \enspace 1 / x \label{MapReflP1LF}
\end{align}

\section{Somos-3: The representation of $T$ in $Cr_2 (\mathbb{R})$ \label{RepSomos3} and invariant curves with its automorphisms}

The terms of the Somos-3 sequence can be expressed in a simple way by the 3 initial values $a_0, a_1, a_3$ as $a_{2n} = a_2^{n} / a_0 ^{n-1}$ and $a_{2n+1} = a_1 a_2^{n} / a_0 ^{n}$ for all $n$. \\  

The transform $a_{n} \mapsto a_{n} \, b c^n$ results in a 1-paramter map depending on $c \ne 0$:
\begin{equation}
M(c) : \enspace  \mathbb{RP}^2 \dashrightarrow \enspace \mathbb{RP}^2, 
 ( a_0 \, : a_1 \, : a_2 ) \mapsto \enspace ( a_0 \, : a_1 c \, : a_2 c^2 ) \label{MapMonomialP2}
\end{equation}

The reflection transform $a_{n} \mapsto a_{2-n}$ results in the reflection map $R$:
\begin{equation}
R : \enspace  \mathbb{RP}^2 \dashrightarrow \enspace \mathbb{RP}^2, 
 ( a_0 \, : a_1 \, : a_2 ) \mapsto \enspace ( a_2 \, : a_1 \, : a_0 ) \label{MapReflP2}
\end{equation}

The shift transform $a_{n} \mapsto a_{n+1}$ results in the shift map $F$:
\begin{equation}
F : \enspace  \mathbb{RP}^2 \dashrightarrow \enspace \mathbb{RP}^2, 
 ( a_0 \, : a_1 \, : a_2 ) \mapsto \enspace ( a_1 \, : a_2 \, : a_1 a_2 / a_0 ) \label{MapShiftP2}
\end{equation}
In this case F is \emph{of finite order}. It has order $2$.

The group $\langle \, R, F \, \rangle$ is therefore isomorph to the finite dihedral group $D_4$.
This group can also be generated
by $R$ and an additional involution $G = R F$, so $\langle \, R, F \, \rangle = \langle \, R, G \, \rangle$: 
\begin{equation}
G = R F : \enspace  \mathbb{RP}^2 \dashrightarrow \enspace \mathbb{RP}^2, 
 ( a_0 \, : a_1 \, : a_2 ) \mapsto \enspace ( a_1 a_2 / a_0 \, : a_2  \, : a_1 ) \label{InvoP2}
\end{equation}
This map stems from the reflection transformation $a_{n} \mapsto a_{3-n}$. This is the standard involution. As polynomial map it is of degree $2$. $\det\mathrm{jac}$ is the union of the 3 lines $a_0 = 0, a_1 = 0, a_2 = 0$. \\

The 1-parameter map $M(c)$ does not commute with the 2 involutions $R, G$. We have $R M(c) R^{-1} = G M(c) G^{-1} = M(c)^{-1} = M(c^{-1})$. So $\langle \, M(c), R \, \rangle$ and $\langle \, M(c), G \, \rangle$ are infinite, \emph{now continuous} dihedral groups $D_\infty (c)$. \\

$\hat{T} = \langle \, M(c), \; R, G \, \rangle \simeq T^{*}_2 \simeq \mathbb{Z}_2 \times \Rx \rtimes \mathbb{Z}_2$, see \ref{GroupSomosTransStar_12}, defines a representation of the group of Somos sequence transformations $T$ in $Cr_2 (\mathbb{R})$.
The degrees of the maps $F^n, R F^n, M(c)$ and all group elements are  $1, 2, 1, 2, 1, 2, 1, 2, \dots$
for $n \ge 0$. The degrees are periodic $mod \; 2$. \\

\subsection{A construction of curves invariant under $\hat{T}$}

An invariant curve of degree $d$ is composed of $n =floor (d/2)+1$ monomials $a_0^m a_1^{d-2m} a_2^m$.
These monomials are fixed by the map $R$. Applying the map $M(c)$ to these monomials, they all aquire the same factor $c^d$ independent form the monomial. Because for odd $d$ all these monomials have the factor $a_1$, all curves obtained by a linear combination the monomials are reducible. \\
So we have to treat only the case of even degree $d$. The map $G$ is the standard involution and maps $a_0^m a_1^{d-2m} a_2^m  \mapsto a_0^{-m} a_1^{2m-d} a_2^{-m}$. Multiplying all images with
$a_0^{2/d} a_1^{d} a_2^{2/d}$ we get the $n$ monomials again. So under $G$ pairs of monomials are permuted,
for $n$ multiples of $4$ one monomial is fixed. Now order the pairs and and and a fixed monomial with decraesing powers of $a_0$, in a pair the $a_0$ with the higher power. A linear combination with a free parameter for each pair except the first (and a fixed monomial) results in an invariant curve. Inserting a $-$ sign in the sum of all pairs of monomials and omitting a fixed monomial we get an antisymmetric (under $G$) version of the curve. \\ 

Example: \\
For the degree $d = 4$ we get the 3 monomials $a_0^0 a_1^4 a_2^0, a_0^1 a_1^2 a_2^1, a_0^2 a_1^0 a_2^2$.
$G$ is permuting $a_1^4$ and $a_0^2 a_2^2$. $G$ fixes $a_0^1 a_1^2 a_2^1$. So a linear combination
of the sum of the permuted pair and the fixed monomial is invariant. This curve has a free parameter.
The under $G$ antisymmetric version of a curve is obtained by the difference of the permuted pair. 

\subsection{Quadratic curves invariant under $\hat{T}$}

These 2 quadratic curves are left invariant by $\hat{T}$:
\begin{equation}
	C_2 = a_{0} a_{2} + a_{1}^2 \label{QuadraticCurveP2}
\end{equation}
\begin{equation}
	C_{2a} = a_{0} a_{2} - a_{1}^2 \label{QuadraticCurveP2a}
\end{equation}

Because $C_2$ and $C_{2a}$ are quadratic in $a_1$ we get another birational involution:
\begin{equation}
H : \enspace  \mathbb{RP}^2 \dashrightarrow \enspace \mathbb{RP}^2, 
 ( a_0 \, : a_1 \, : a_2 ) \mapsto \enspace ( a_0 \, : - a_1 \, : a_2 ) \label{MapExtraP2}
\end{equation}
This map as a linear map is of degree $1$ and commutes with the 1-parameter map $M(c)$.
Because $C_2^2 = C_4 (+2)$ and $C_{2a}^2 = C_4 (-2)$, the map $J (\pm 2)$, see \ref{MapExtraP2}
leaves $C_{2}$ and $C_{2a}$ invariant too. \\ 

The group $\langle \, \hat{T}, H, \, J (\pm 2) \, \rangle$ of automorphisms of $C_2$ and $C_{2a}$ is
isomorph to \\ $\mathbb{Z}_2^2 \times (\Rx \times \mathbb{Z}) \rtimes \mathbb{Z}_2$. The degrees of the maps
in this group are $1, 2, 4, 6, 8, \dots$.

\subsection{Quartic curves invariant under $\hat{T}$}

This quartic curve is left invariant by $\hat{T}$:
\begin{equation}
	C_4 (\alpha) = (a_{0}^2 a_{2}^2 + a_{1}^4) + \alpha \, a_0  a_{1}^2 a_2 \label{QuarticCurve}
\end{equation}
For $\alpha = -2, +2$ this curve $C_4$ is reducible: $C_4 (+2) = C_2^2, C_4 (-2) = C_{2a}^2$.   
For $\alpha = 0$ this curve $C_4$ is reducible in $\mathbb{C}$: $C_4 (0) =  (a_{0} a_{2} + i a_{1}^2) \, (a_{0} a_{2} - i a_{1}^2)$. Another anti-symmetric invariant curve is $C_{4a} = a_{0}^2 a_{2}^2 - a_{1}^4$ which is reducible in $\mathbb{R}$ as $C_{4a} = C_2 \, C_{2a}$. \\

Because $C_4 (\alpha)$ contains $a_1$ only with even powers we get another involution:
\begin{equation}
H_1 : \enspace  \mathbb{RP}^2 \dashrightarrow \enspace \mathbb{RP}^2, 
 ( a_0 \, : a_1 \, : a_2 ) \mapsto \enspace ( a_0 \, : - a_1 \, : a_2 ) \label{MapExtraP2}
\end{equation}
Because $C_4 (\alpha)$ is quadratic in $a_0$ and $a_2$ we can construct two further birational involutions.
\begin{equation}
H (\alpha) : \enspace  \mathbb{RP}^2 \dashrightarrow \enspace \mathbb{RP}^2, 
 ( a_0 \, : a_1 \, : a_2 ) \mapsto \enspace ( - \alpha \, a_1^2 / a_2 - a_0 \, : a_1 \, : a_2 ) \label{MapExtraP2}
\end{equation}
This map as polynomial map is of degree $2$ and commutes with the 1-parameter map $M(c)$.
$\det\mathrm{jac}$ is the line $a_2 = 0$ with multiplicity $3$. The other involution belonging to $a_2$ is just this $H$ conjugated by $R$. \\

The group $\langle \, \hat{T}, H_1, H (\alpha) \, \rangle$ of automorphisms of $C_4$ is $\mathbb{Z}_2^2 \times (\Rx \times \mathbb{Z}) \rtimes \mathbb{Z}_2$. Here the group $\mathbb{Z}$ in this direct product is generated by $J (\alpha) = R H (\alpha)$.

With $U (\alpha) = a_0 a_2 + \alpha \, a_1^2$ we get 2 forms for all even and odd powers of $J$: \\ 
\begin{align}
J^{2n} (\alpha) & : \dots , 
 a_{0123} \mapsto \enspace ( U^{2n}(\alpha) \, & : (-1)^n a_0^{n-1} a_1 a_2^{n} \, U^{n}(\alpha) \, & : a_0^{2n-1} a_2^{2n+1} ) \label{MapExtraP2} \\
J^{2n+1} (\alpha) & : \dots , 
 a_{0123}  \mapsto \enspace ( a_0^{2n} a_2^{2n+2} \, & : (-1)^n a_0^{n} a_1 a_2^{n+1} \, U^{n}(\alpha) \, & : - U^{2n+1}(\alpha) ) \label{MapExtraP2}
\end{align}

The degrees of the maps in  $\langle \, \hat{T}, H_1, H (\alpha) \, \rangle$ are $1, 2, 4, 6, 8, \dots$.

\subsection{Sectic curves invariant under $\hat{T}$}

In a similar manner as the 2 quartic curves were constructed we get 2 sectic curves left invariant by $\hat{T}$. But the 2 curves are reducible. $C_6 (\alpha) = (a_{0}^3 a_{2}^3 + a_{1}^6) + \alpha \, (a_{0}^2 a_{1}^2 a_{2}^2 + a_{0}^1 a_{1}^4 a_{2}^1)$ factorizes as $C_6 (\alpha) = C_4 (\alpha - 1) \, C_2$. the
antisymmetric $C_{6a} (\alpha) = (a_{0}^3 a_{2}^3 - a_{1}^6) + \alpha \, (a_{0}^2 a_{1}^2 a_{2}^2 - a_{0}^1 a_{1}^4 a_{2}^1 )$ factorizes as $C_{6a} (\alpha) = C_4 (\alpha - 1) \, C_{2a}$.   

\subsection{Octic curves invariant under $\hat{T}$}

In a similar manner as the quartic curve is constructed we get these octic curves left invariant by $\hat{T}$:
\begin{equation}
	C_8 (\alpha, \beta) = (a_{0}^4 a_{2}^4 + a_{1}^8) + \alpha \, (a_{0}^3 a_{1}^2 a_{2}^3 + a_{0}^1 a_{1}^6 a_{2}^1) + \beta \, a_{0}^2 a_{1}^4 a_{2}^2\label{OcticCurve}
\end{equation}
This curve is reducible in $\mathbb{C}$. For $\alpha^2 - 4 (\beta - 2) \ge 0$ this curve $C_8$ is already reducible in $\mathbb{R}$. \\
This because $C_8 (\alpha, \beta) = C_8 (x_1 + x_2, \, 2 + x_1  x_2) = C_4 (x_1) \, C_4 (x_2)$, we can determine $x_1, x_2$ solving the quadratic equation $x^2 - \alpha \, x + (\beta - 2) = 0$ and so we get a factorization in $\mathbb{R}$ or $\mathbb{C}$ depending on the sign of the discriminant.
\begin{equation}
	C_{6a} (\alpha) = (a_{0}^3 a_{2}^3 - a_{1}^6) + \alpha \, (a_{0}^2 a_{1}^2 a_{2}^2 - a_{0}^1 a_{1}^4 a_{2}^1 ) \label{OcticCurve}
\end{equation}
This antisymmetric octic is reducible: $C_{8a} (\alpha) = C_4 (\alpha) \, C_{2} \, C_{2a}$.  

\subsection{Curves with degree $> 8$ invariant under $\hat{T}$}

$C_{4k} (\dots)$ and $C_{4k+2} (\dots)$ have $k$ arguments. we have the following factorization:
$$C_{4k} (\alpha_1, \dots, \alpha_k) = \prod_{i=1}^{k} C_4 (x_i)$$
$$C_{4k+2} (\alpha_1, \dots, \alpha_k) = C_2 \, C_{4k} (\alpha_1, \dots, \alpha_k)$$
The $\alpha_1, \dots, \alpha_k$ are sums of elementary symmetric functions of $x_1, \dots, x_k$
So $C_{4k}$ is reducible in $\mathbb{C}$, if the polynomial in $\alpha_1, \dots, \alpha_k$ has a real root
it is already reducible in $\mathbb{R}$. $C_{4k+2}$ is reducible because a polynomial odd degree has a real root. \\
An example: With these elementary symmetric functions $\sigma_1 = x_1 + x_2 + x_3, \sigma_2 = x_1 x_2 + x_1 x_3 + x_2 x_3, \sigma_3 = x_1 x_2 x_3$, we get $C_{12} (\alpha_1, \alpha_2, \alpha_3) = C_{12} (\sigma_1, 3 + \sigma_2, 2 \sigma_1 + \sigma_3) = C_4 (x_1) \, C_4 (x_2) \, C_4 (x_3)$. The corresponding polynomial is
$x^3 - \alpha_1 \, x^2 + (\alpha_2 - 3) \, x - (\alpha_3 - 2 \alpha_1) = 0$.  

\section{Somos-5: The representation of $T$ in $Cr_4 (\mathbb{R})$ \label{RepSomos5}}

The transform $a_{n} \mapsto a_{n} \, b c^n$ results in a 1-paramter map depending on $c \ne 0$:
\begin{equation}
M(c) : \enspace  \mathbb{RP}^4 \dashrightarrow \enspace \mathbb{RP}^4, 
 ( a_0 \, : a_1 \, : a_2 \, : a_3 \, : a_4 ) \mapsto \enspace ( a_0 \, : a_1 c \, : a_2 c^2  \, : a_3 c^3  \, : a_4 c^4 ) \label{MapMonomialP4}
\end{equation}

The reflection transform $a_{n} \mapsto a_{4-n}$ results in the reflection map $R$:
\begin{equation}
R : \enspace  \mathbb{RP}^4 \dashrightarrow \enspace \mathbb{RP}^4, 
 ( a_0 \, : a_1 \, : a_2 \, : a_3 \, : a_4 ) \mapsto \enspace ( a_4 \, : a_3 \, : a_2 \, : a_1 \, : a_0 ) \label{MapReflP4}
\end{equation}

The shift transform $a_{n} \mapsto a_{n+1}$ results in the shift map $F$:
\begin{equation}
F : \enspace  \mathbb{RP}^4 \dashrightarrow \enspace \mathbb{RP}^4, 
 ( a_0 \, : a_1 \, : a_2 \, : a_3 \, : a_4 ) \mapsto \enspace ( a_1 \, : a_2 \, : a_3 \, : a_4  \, : (a_1 a_4 + a_2 a_3) / a_0 ) \label{MapShiftP4}
\end{equation}

The group $\langle \, R, F \, \rangle$ is isomorph to the infinite dihedral group $D_\infty$.
This group can also be generated
by $R$ and an additional involution $G = R F$, so $\langle \, R, F \, \rangle = \langle \, R, G \, \rangle$: 
\begin{equation}
G = R F : \enspace  \mathbb{RP}^4 \dashrightarrow \enspace \mathbb{RP}^4, 
 a_{01234} \mapsto \enspace ( (a_1 a_4 + a_2 a_3) / a_0 \, : a_4 \, : a_3 \, : a_2  \, : a_1 ) \label{InvoP4}
\end{equation}
This map stems from the reflection transformation $a_{n} \mapsto a_{5-n}$. As polynomial map it is of degree $2$. $\det\mathrm{jac}$ is the union of the hyperplane $a_0 = 0$ and the quadratic hypersurface $a_1 a_4 + a_2 a_3 = 0$. \\

The 1-parameter map $M(c)$ does not commute with the 2 involutions $R, G$. We have $R M(c) R^{-1} = G M(c) G^{-1} = M(c)^{-1} = M(c^{-1})$. So $\langle \, M(c), R \, \rangle$ and $\langle \, M(c), G \, \rangle$ are infinite, \emph{now continuous} dihedral groups $D_\infty (c)$. \\

$\hat{T} = \langle \, M(c), \; R, G \, \rangle \simeq T^{*} \simeq (\Rx \times \mathbb{Z}_2) \rtimes \mathbb{Z}_2$, see \ref{GroupSomosTransStar}, defines a representation of the group of Somos sequence transformations $T$ in $Cr_4 (\mathbb{R})$.
The degrees of the maps $F^n, R F^n, M(c)$ and all group elements are $1, 2, 3, 4, 6, 9, 11, \dots$
for $n \ge 0$.  \\

Question: Are there 3-dimensional hypersurfaces invariant under $\langle \, M(c), R, G  \, \rangle$ or a nontrivial subgroup? \\

\bibliographystyle{amsplain}

\begin{thebibliography}{10}
\bibitem{Deserti} J.~D\'{e}serti: \emph{Some Properties of the Cremona Group}. Ensaios Matem{\'a}ticos \textbf{21} (2012) pp. 1--188, \href {https://doi.org/10.21711/217504322012/em211} {\path{DOI: 10.21711/217504322012/em211}}.
\bibitem{EPSW} G. Everest, A. van der Poorten, I. Shparlinski and T. Ward: \emph{Recurrence Sequences}.
 Amer. Math. Soc. (2003) pp. 9--179.
\bibitem{FoZe} S. Fomin and A. Zelevinsky: \emph{The Laurent phenomemon}. arXiv:math/0104241 [math.CO], 2001,
\href{https://doi.org/10.48550/arXiv.math/0104241}{DOI: 10.48550/arXiv.math/0104241}
\bibitem{Gale} D. Gale: \emph{The strange and surprising saga of the Somos sequences}. Mathematical Intelligencer \textbf{13} (1991) pp. 40--43, \href {https://doi.org/10.1007/BF03024070} {\path{DOI: 10.1007/BF03024070}}.
\bibitem{Sloane} N.~J.~A. Sloane: \emph{On-Line Encyclopedia of Integer Sequences}.
 \href {https://oeis.org} {\path{Website}}.
\end{thebibliography}

\end{document}